\RequirePackage{fix-cm}
\documentclass[smallextended,preprint]{svjour3}       
\smartqed  
\usepackage{graphicx}
\usepackage{amsmath}
\usepackage{amssymb}
\usepackage{amsfonts}
\usepackage{graphicx}
\usepackage[utf8]{inputenc}
\usepackage{enumerate}
\usepackage{enumitem}
\usepackage[english]{babel}
\usepackage{bm}
\usepackage{pdfpages}
\usepackage{float} 
\usepackage{latexsym}
\usepackage{verbatim} 
\usepackage{color}
\usepackage{bbm}
\usepackage{setspace}
\usepackage{cancel}
\usepackage{units}
\usepackage{multirow}
\usepackage{fancyhdr}
\usepackage[numbers]{natbib}
\usepackage[12pt]{extsizes}
\usepackage{geometry}
\newcommand{\E}[1]{\mathbb{E}\left[#1\right]}
\newcommand{\x}{\textbf{x}}
\newcommand{\n}{\textbf{n}}
\newcommand{\m}{\textbf{m}}
\newcommand{\e}{\textbf{e}}
\newcommand{\N}{\textbf{N}}
\newcommand{\X}{\textbf{X}}
\newcommand{\h}{\textbf{h}}
\newcommand{\ud}{\textrm{d}}
\newcommand{\trans}{\mathrm{T}}
\newcommand{\norm}[1]{\left\lVert#1\right\rVert_1}

\journalname{Journal of Mathematical Biology}
\begin{document}


\title{A dual process for the coupled Wright-Fisher diffusion
}


\author{Martina Favero        \and
        Henrik Hult           \and
        Timo Koski
}

\authorrunning{M.Favero, H.Hult, T.Koski} 

\institute{  {M. Favero }
	         \at
             Department of Mathematics, KTH, 11428 Stockholm, Sweden \\
             \email{mfavero@kth.se}           
	         \and
			 H. Hult \at
	         Department of Mathematics, KTH, 11428 Stockholm, Sweden \\ 
	         \email{hult@kth.se}  
			 \and 
	         T. Koski \at
		     Department of Mathematics, KTH, 11428 Stockholm, Sweden \\ 
	         \email{tjtkoski@kth.se}  
}

\date{  }

\maketitle

\begin{abstract}
The coupled Wright-Fisher diffusion is a multi-dimensional Wright-Fisher diffusion for multi-locus and multi-allelic genetic frequencies, expressed as the strong solution to a system of stochastic differential equations that are coupled in the drift, where the pairwise interaction among loci is modelled by an inter-locus selection. In this paper, an ancestral process, which is dual  to the coupled Wright-Fisher diffusion, is derived. The dual process corresponds to the block counting process of coupled ancestral selection graphs, one for each locus. Jumps of the dual process arise from coalescence, mutation, single-branching, which occur at one locus at the time, and double-branching, which occur simultaneously at two loci. The coalescence and mutation rates have  the typical structure of the transition rates of the Kingman coalescent process. The single-branching rate  not only contains the one-locus selection parameters  in a form that generalises the rates of an ancestral selection graph, but it also contains the two-locus selection parameters to include the effect of the pairwise interaction on the single loci. The double-branching rate reflects the particular structure of pairwise selection interactions of the coupled Wright-Fisher diffusion. 
Moreover, in the special case of two loci, two alleles, with selection and parent independent mutation, the stationary density  for the coupled Wright-Fisher diffusion  and  the transition rates of the dual process are obtained in an explicit form.


\end{abstract}

\section{Introduction}

The coupled Wright-Fisher diffusion  was introduced by Aurell, Ekeberg and Koski \cite{timo} with the purpose of analysing networks of loci in recombining populations of bacteria,  or more precisely,  detecting couples of loci co-evolving under strong selective pressure when the linkage disequilibrium is low across the genome.  The model includes parent dependent mutation, interlocus selection and free recombination. Mutation is assumed to occur independently at each locus, while selection consists of first and second order selective interaction among loci.

This particular type of assumptions on the selection and recombination structure are suitable for example for some populations of bacteria, as showed in \cite{Skwark2017},
where the same type of assumptions are made.
In \cite{Skwark2017}, it is explained that the high amount of homologous recombination in populations of \textit{Streptococcus Pneumoniae}, which results in low linkage disequilibrium across the genome, makes this population ideal for detecting genes that evolve under shared selection pressure.
On the contrary, in other populations of bacteria, e.g. \textit{Streptococcus Pyogenes}, the low amount of homologous recombination makes it difficult to separate couplings attributable to recombination from those attributable to selection and thus the assumptions above are not suitable to study such populations. 

The mathematical idea corresponding to these biological characteristics, is that the recombination is high enough to be approximated with infinite recombination, which would make the processes at each locus independent,  and it is thus the selection only that causes the coupling between the diffusions at the different loci.

Furthermore, it is assumed that selection acts on the individual loci and on pairs of loci.
The pairwise selection can be thought of as a network, where the vertices represent the loci and the edges the possible interactions, as shown in \cite{timo}.
Of course, the possible set of interactions could, in principle, be more complex than a network, but considering pairwise interactions turns out to be useful to reveal certain types of co-evolutionary patterns, see \cite{Skwark2017}.

The model considers $L$ different loci where,  at each locus,  a number of  variants (alleles) are possible.  The allele types at locus $l$ are labelled by $1,\dots,M_l$, thus assuming that the type space at each locus is finite. The population is haploid. The coupled Wright-Fisher diffusion is obtained  as the weak limit of a sequence of discrete Wright-Fisher models characterised by the assumption that the evolution of the population at one locus is conditionally independent of the other loci given that the previous generation at each locus is known, see \cite{timo} for details. 
It is based on quasi-linkage equilibrium where the fitness coefficients, see Section \ref{notation}, are inspired by a Potts model, see \cite{Aurell2019,Shraiman2011}, and generalise the classical additive fitness under weak selection, see e.g. \cite[Ch. II]{Burger2000}, to the multi-locus case.  With two loci and without the fist order selection terms, the coupled Wright-Fisher diffusion is reduced to a haploid version of the model with weak selection, loose linkage in \cite{ethier1989}.

Here we state the definition of the diffusion as the solution of a system of stochastic differential equations, without reference to the underlying discrete model.
The coupled Wright-Fisher diffusion, $\mathbf{X} = \{\mathbf{X}(t), t\geq 0\}$,  represents the evolution of the vector of all frequencies  of allele types at each locus. Let 
$$\mathbf{X}^{(l)}(t)=(X_1^{(l)}(t),\dots,X_{M_l}^{(l)}(t))^\trans$$ 
represent the vector of frequencies at locus $l$, with $X_i^{(l)}(t)$ being the frequency of allele type $i$ at locus $l$, then $$\mathbf{X}(t)=(\mathbf{X}^{(1)}(t)^\trans,\dots,\mathbf{X}^{(L)}(t)^\trans)^\trans.$$

The process $\mathbf{X} $ is the strong solution to the system of stochastic differential equations
\begin{equation}
\label{SDE}
\ud \mathbf{X}(t)= 
\mu(\mathbf{X}(t)) \ud t +
D(\mathbf{X}(t)) \nabla V(\mathbf{X}(t)) \ud t +
D^{\nicefrac{1}{2}}(\mathbf{X}(t)) \ud \mathbf{W}(t),
\end{equation}
where $V$ is a specific quadratic function encoding the structure of the interactions, $\nabla V$ its gradient, 
while the mutation  vector $\mu$ and the diffusion matrix  $D$ have the following block structure,
\begin{align*}
\mu(\mathbf{x})=
\left(
\begin{matrix}
\mu^{(1)}(\x^{(1)})\\
\vdots \\
\mu^{(L)}(\x^{(L)})
\end{matrix}
\right),
\quad\quad
D(\mathbf{x})=
\left(
\begin{matrix}
D^{(1)}(\x^{(1)})&  		 &         \\
& \ddots 	 &         \\
&           & D^{(L)}(\x^{(L)})
\end{matrix}
\right),
\end{align*}
with
$\mu^{(l)}:\mathbb{R}^{M_l}\rightarrow \mathbb{R}^{M_l}$ and 
$D^{(l)}:\mathbb{R}^{M_l}\rightarrow \mathbb{R}^{M_l \times M_l}$.
The functions 
$V$, $\mu$ and  $D$ 
are described in detail in the next section. The process $\mathbf{W}=((\mathbf{W}^{(1)})^\trans ,\dots,(\mathbf{W}^{(L)})^\trans)^\trans$ is a multidimensional Brownian motion with $\mathbf{W}^{(l)}$ having the dimension of $\mathbf{X}^{(l)}$.

The system of SDEs (\ref{SDE}) consists of
$L$ systems of equations for  $\mathbf{X}^{(1)},\dots,\mathbf{X}^{(L)}$, coupled by the drift term  $D\,\nabla V$.
Note that, if $\nabla V = 0$, there is no interaction among the loci and the coupled Wright-Fisher diffusion consists of $L$ independent Wright-Fisher diffusions, that is,  each $\mathbf{X}^{(l)}$ solves
$$
\ud \mathbf{X}^{(l)}(t)= 
\mu^{(l)}(\mathbf{X}^{(l)}(t)) \ud t +
D^{(l)\nicefrac{1}{2}}(\mathbf{X}^{(l)}(t))\, \ud \mathbf{W}^{(l)}(t),
$$
which is the SDE for a single-locus, multi-type Wright-Fisher diffusion with mutations. 
In fact, the coupling of the loci is entirely due to selective interactions that are described by the drift term $D\,\nabla V$. Without the interaction drift term, the diffusion in this paper, with $L=2$, is reduced to the independent-loci model in \cite{ethier1990}.  That is,  the weak  limit of a sequence of multi-locus neutral Wright-Fisher diffusions with recombination rate going to infinity. In the multi-locus case, the same diffusion  appears also in \cite[Sect. 3.3]{Griffiths2016} as an example under free recombination.

An interesting feature of the coupled Wright-Fisher diffusion, addressed by Aurell et al. \cite{timo} as one of the main motivations for its introduction, is its stationary density which appeared, in a more general
form, as a conjecture by Kimura over half a century ago. In \cite{Kimura1955}, Kimura suggests a Wright-Fisher model for multi-locus and multi-allelic genetic frequencies and conjectures that the stationary density is of the form 
$\pi e^m$,
where $\pi$ is the product of Dirichlet densities and $ m $ is a generic mean fitness term. 
The coupled Wright-Fisher diffusion is constructed so that the quadratic function $ V $ could replace the generic $ m $. Indeed, under the assumption of parent independent mutations, the stationary
density, $p$, of the coupled Wright-Fisher diffusion is known up to a normalising constant $Z$, and corresponds to  the one conjectured by Kimura with $ m = 2V $,
\begin{equation}
\label{stationarydensity}
p \propto \pi  e^{2 V}, 
\end{equation}
see Section \ref{notation} for the definition of $\pi$ and $V$.  In fact, the form of the stationary density, under parent independent mutations, relies on the fact that the covariance of the diffusion defines a Svirezhev-Shahshahani metric on the simplex, with respect to which the drift is a gradient, see \cite[Appendix E.3]{Burger2000}. 

In this paper a dual process for the coupled Wright-Fisher diffusion is studied. In population genetics, Markov duality has proven its effectiveness in combining information from two processes related to the same population: a diffusion process modelling the evolution of frequencies of genetic types forward in time and a reverse-time jump process modelling the ancestral history of a sample of individuals taken at the present time. The simplest and most well known duality relationship in this context is the moment duality between the Wright-Fisher diffusion and the block counting process of the Kingman coalescent.

The strength  of Markov duality is that it provides a tool to analyse properties of the population by combining knowledge about the forward-in-time process and the backward-in-time process.
Even  when both processes are complicated,  
as often happens when mutation, recombination or selection mechanisms are involved, some known properties of one process can be used to analyse unknown properties of the other process and vice versa, leading to further insights about the population. 

Several duality relationships have been established between various generalisations of the Wright-Fisher diffusion and the associated time reversed ancestral processes generalising the coalescent process.
For example, when the selection mechanism is taken into account, the ancestral process associated to the Wright-Fisher diffusion is the ancestral selection graph (ASG), see \cite{Krone1997, Neuhauser1997}, which is closely related to the dual process in this paper when only one locus is considered, see Section \ref{ASG}. Unlike the Kingman coalescent, which has a tree structure, the ASG is branching and coalescing: the ancestral tree is replaced by an ancestral graph containing true and virtual lineages and embedding the genealogy of the sample of individuals. For a complete survey on duality for Markov processes, see \citep{Jensen2014}, and for a brief  overview  of duality in population genetics  see \cite{Griffiths2016} and the references therein.

In this paper, the main result concerns the derivation of a dual process for the coupled Wright-Fisher diffusion. The results show that, in this model, the dual process corresponds to the block counting process of $L$ coupled ASGs, one for each locus, evolving simultaneously. Coalescence, mutation and single-branching, which is due to selection acting on the single loci, occur at different times in the different ASGs, whereas branching that is due to selection acting on pairs of loci, occurs simultaneously in two ASGs. The latter type of branching is referred to as double-branching in this paper. The main result in this paper is Theorem \ref{thm:main}, which provides a description of the transition rates of the pure jump Markov process, $\mathbf{N} = \{\mathbf{N}(t), t\geq 0\}$, that is dual to the coupled Wright-Fisher diffusion, $\mathbf{X} $, through the  duality relationship
\begin{equation}
\label{duality_expectations}
\E{F(\mathbf{X}(t),\n)|\mathbf{X}(0)=\x}=
\E{F(\x,\mathbf{N}(t))|\mathbf{N}(0)=\n},
\end{equation}
where $F$ is a duality function, to be determined. The derivation uses a generator approach as in \cite{Griffiths2016} and \citep{Etheridge2009}. It is based on the duality relationship of the infinitesimal generators
\begin{equation}
\label{duality_generators}
\mathcal{L}F(\cdot,\n)(\x)=\mathcal{L}^D F(\x,\cdot \,)(\n),
\end{equation}
where $\mathcal{L}$ is the generator of the coupled Wright-Fisher diffusion and $\mathcal{L}^D$ the unknown generator of the dual process. 
By proposing an appropriate duality function $F$, the generator $\mathcal{L}^D$ of the dual process can be identified, from which transition rates of the dual process are obtained. 
Under mild conditions, which are verified in Section \ref{proof},
the method of duality \cite[Ch. 4]{Ethier1986}, also used in e.g. \cite{barbour2000,Etheridge2009,Mano2009}, ensures that
 the duality relationship of the generators (\ref{duality_generators}) implies (\ref{duality_expectations}).

Understanding the structure of an ancestral process, $\N$, which is dual to  a diffusion of the Wright-Fisher type,  $\X$, plays a significant role in population genetics inference. 
As is often the case, the available data consist of  observations of the genetic types of a sample of individuals at the present time, $\N(0)=\n$, whereas the evolution of the process is not observed. This results in the likelihood function being intractable when the size of the population is large.
In order to compute the likelihood, one could, in theory, condition on the genealogical history of the sample and then integrate over all possible histories that are compatible with the sample. However, the domain of integration is so large that, in practice, numerical integration methods are useless even for intermediate sized populations. Simulation-based methods are generally preferred.  As carefully explained by Stephens in \cite{Stephens2007},
naive Monte Carlo methods based on simulating the histories forward in time produce next to useless approximations of the likelihood for problems involving samples of a more than few individuals. This is due to the fact that only very few simulations contribute significantly to the approximation, 
while the contribution of the remaining simulations is negligible.
Simulation- and likelihood-based techniques that have proven to work for these problems are  Markov chain Monte Carlo,  importance sampling and sequential Monte Carlo. All these methods rely on knowing, to some extent, the structure of the ancestral process in order to approximate its backward dynamics, see e.g.\  \cite{Stephens2007, griffiths1994, Stephens2000, stephens2003, Koskela2015, Koskela2018} for details. 

From the duality relation (\ref{duality_expectations}), it is also possible to derive an expansion of the transition distribution of the diffusion $\X$, see \cite{Etheridge2009,Griffiths2016,barbour2000}, in terms of the limit of the transition densities of the dual process $\N$. 
In the absence of mutation, the duality relation \eqref{duality_expectations} can also be used to determine fixation probabilities. That is, the probability that the frequency of a given allele at a given loci is equal to $1$. Such probabilities may be studied by taking the limit, as $t \to \infty$, in \eqref{duality_expectations} and considering the recurrence/transience properties of the dual process $\mathbf{N}$, see e.g.\ \cite{Mano2009, Foucart2013, Griffiths2016, gonzalez2018} for studies of  the Wright-Fisher process with selection and frequency dependent selection,  the multi-locus Wright-Fisher process with recombination, the $\Lambda$-Wright-Fisher process with selection and the $\Xi$-Wright-Fisher process with frequency dependent selection, respectively. 

The paper is outlined as follows. 
In Section \ref{notation} a background on  the coupled Wright-Fisher diffusion is provided. 
Section \ref{outline} outlines the  general generator approach to derive a dual process. 
In Section \ref{ASG} the case of one locus, two allele types and  parent independent mutations is considered. In this case the dual process is related to the ancestral selection 
graph, moreover,  explicit formulas for the stationary density of the diffusion and the transition rates of the dual process are obtained. 
The main result is provided in Section \ref{multi}, and proved in Section \ref{proof}, where a dual process is derived in the general multi-locus setting. The final Section 
\ref{ex} provides additional details in the case of two loci, two alleles, selection and parent independent mutations, more precisely, the transition rates of the dual process are expressed in terms of beta and confluent hypergeometric functions.

\section{Preliminaries on the coupled Wright-Fisher diffusion} 
\label{notation}

In this section the coupled Wright-Fisher diffusion is introduced and the explicit expression for its infinitesimal generator is provided.  The notation in this section differs slightly from that in \cite{timo}, where 
the frequency of the last allele type at each locus is omitted, being a function of the other frequencies, whereas in this paper an expanded version of the diffusion is considered, which includes all the frequencies. Since the frequencies sum up to one the descriptions are equivalent. For our purpose we find the expanded version more convenient to work with.

For a given integer $L\geq 1$, the number of loci, let $M_1, \dots, M_L$ be positive integers representing the number of alleles at each locus. Put $M=\sum_{l=1}^{L}M_l$.
A vector $\x\in [0,1]^{M}$ is interpreted as the concatenation of 
$L$ vectors with lengths $M_1,\dots,M_L$, i.e. $\x = ((\x^{(1)})^\trans, \dots, (\x^{(L)})^\trans)^\trans$ with $\x^{(l)}\in [0,1]^{M_l}$, $ l=1,\dots L $, and the coordinate $i$ in vector $\x^{(l)}$ is denoted by $x _{i}^{(l)}$.
Similarly, a matrix $A\in \mathbb{R}^{M\times M}$ consists of $L^2$
blocks with dimensions $(M_l\times M_r)_{l,r=1,\dots,L}$. The block at position $(l,r)$ is denoted by 
$A^{(lr)}$ and its  component at position $(i,j)$ is denoted by by $A_{ij}^{(lr)}$.
Furthermore   $\mathbf{e}_{i}^{(l)}$  denotes the unit vector in $\mathbb{R}^M$ with the $i^{th}$ component of its $l^{th}$ building vector being equal to $1$.

In the following, each of the terms appearing  in (\ref{SDE}) will be described, starting from the interaction drift term.
The quadratic function $V:[0,1]^M\rightarrow \mathbb{R}$  is given by
$$
V(\mathbf{x})=
\mathbf{x}^\trans\mathbf{h}+\frac{1}{2}\mathbf{x}^\trans J \mathbf{x}, 
$$
where $\mathbf{h}\in \mathbb{R}^{M}_+$ and
$J \in \mathbb{R }^{M\times M}_+$ is a symmetric block matrix with 
the blocks on  the diagonal equal to zero matrices, i.e. $J^{(ll)} =0 \in \mathbb{R}^{M_l\times M_l}$ and $J^{(lr)}=(J^{(rl)})^\trans$ for all $l,r=1,\dots,L$.
The vector $\mathbf{h}$ and matrix $J$ contain the selection parameters, expressing,   respectively,  the one-locus selection and the selective interaction among pairs of  loci.
In order to clarify the role of the selection parameters in terms of fitness, we may express the fitness coefficient of the haplotype $\sigma=(i_1,\dots,i_L)$ as
    \begin{equation}
    \label{fitness}
    w_\sigma= 1+ \sum_{l=1}^L h_{i_l}^{(l)} +
    \sum_{l=1}^L\sum_{\substack{r=1\\l<r}}^L J_{i_l i_r}^{(rl)}.
    \end{equation}
    
Note that $\nabla V(\mathbf{x})=\mathbf{h}+J\mathbf{x}$,
since the matrix $J$ is symmetric. Let $g(\mathbf{x})=D(\mathbf{x})\nabla V(\mathbf{x})$. Then,  the components of $g(\mathbf{x})$ are
\begin{equation}
\label{g}
g_i^{(l)}(\mathbf{x})=
\sum_{k=1}^{M_l}
d_{ik}^{(l)}(\x^{(l)})\tilde{h}_k^{(l)}(\mathbf{x}),
\quad
\text{with}
\quad
\tilde{h}_k^{(l)}(\mathbf{x})=h_k^{(l)}+
\sum_{\substack{r=1\\r\neq l}}^{L}\sum_{m=1}^{M_r}
J_{km}^{(l r)}x_m^{(r)}.
\end{equation}
The drift function $\mu$ models the mutations. 
It is assumed that mutations occur independently at each locus, in particular, at the $l^{th}$ locus the mutation rate is $\frac{\theta_l}{2}$ and the probability matrix of mutations is
$P^{(l)}=(P_{ij}^{(l)})_{ i,j=1\dots,M_l}$.
The transition rates of mutations from type $i$ to type $j$ at locus $l$ are thus  
$u_{ij}^{(l)}=\frac{\theta_l}{2}P_{ij}^{(l)}$.
As in the standard Wright-Fisher model with parent dependent mutations,
the components of the drift function are defined by
\begin{equation}
\label{mu}
\mu_i^{(l)}(\x^{(l)})=
\sum_{j=1}^{M_l}
[u_{ji}^{(l)}x_j^{(l)}-u_{ij}^{(l)}x_i^{(l)}].
\end{equation}
Finally, the components of the diagonal block $D^{(l)}(\x^{(l)})$ of the diffusion matrix $D(\x)$ are defined by
\begin{equation}
\label{d}
d_{ij}^{(l)}(\x^{(l)})=
x_i^{(l)}(\delta_{ij}-x_j^{(l)})
\quad\text{with}\quad
\delta_{ij}=
\begin{cases}
1 \text{ if } i=j, \\
0 \text{ if } i\neq j, 
\end{cases}
\end{equation}
which is  characteristic for  Wright-Fisher processes.

Having defined $\mu, D, $ and $V$, a compact definition of the  coupled Wright-Fisher diffusion can be given, in terms of its infinitesimal generator. The coupled Wright-Fisher diffusion $\{\mathbf{X}(t)\}_{t\geq 0}$ is a $M$-dimensional diffusion process on the state space 
$$
\mathcal{S}=
\left\{
\mathbf{x}\in [0,1]^M 
\text{ s.t. }
\sum_{i=1}^{M_l}x_i^{(l)}=1, \quad \forall l=1,\dots,L, 
\right\}, 
$$ 
with generator
\begin{equation}
\label{generator}
\begin{aligned}
&\mathcal{L}f(\x)=
\sum_{l=1}^{L}\left[
\sum_{i=1}^{M_l}
\left(\mu_i^{(l)}(\x^{(l)})+g_i^{(l)}(\x)\right)
\frac{\partial f}{\partial x_i^{(l)} }(\x)
+\frac{1}{2}\sum_{i,j=1}^{M_l}
d_{ij}^{(l)}(\x^{(l)})
\frac{\partial ^2 f}{\partial x_i^{(l)} \partial x_j^{(l)} }(\x)
\right], 
\end{aligned}
\end{equation}
where $\mu$, $g$ and $d$ are given by (\ref{mu}), (\ref{g}) and (\ref{d}), respectively.
The generator $\mathcal{L}$ is defined on the domain
$ C^2(\mathcal{S})$.
\

Before proceeding with the derivation of the dual process, the stationary density (\ref{stationarydensity}) is explicitly presented.
Consider representing the coupled Wright-Fisher diffusion on the state space
$$
\mathcal{\bar{S}}=
\left\{
\mathbf{\bar{x}}\in [0,1]^{M-L}
\text{ s.t. }
\sum_{i=1}^{M_l-1}\bar{x}^{(l)}_i \leq 1 \quad \forall l=1,\dots,L
\right\}, 
$$
where 
\begin{align*}
\bar \x = ({\bar x}^{(1)}_1, \dots, {\bar x}^{(1)}_{M_1-1}, \dots, {\bar x}^{(L)}_1, \dots, {\bar x}^{(L)}_{M_L-1})^\trans
\in \mathcal{\bar S}
\end{align*}
is identified with
\begin{align*}
\x = \left({\bar x}^{(1)}_1, \dots, {\bar x}^{(1)}_{M_1-1}, 1-\sum_{i=1}^{M_1 -1} {\bar x}^{(1)}_i,  \dots, {\bar x}^{(L)}_1, 
\dots, {\bar x}^{(L)}_{M_L-1}, 1-\sum_{i=1}^{M_L -1} {\bar x}^{(L)}_i\right)^\trans \in \mathcal{S}.
\end{align*}
If there are no interactions among loci,  the coupled Wright-Fisher diffusion consists of $L$ independent Wright-Fisher diffusions and the stationary density is well known when the mutations are parent independent.  Wright himself  proved that the stationary distribution of a single-locus, multi-type  Wright-Fisher diffusion with parent independent mutations is Dirichlet, see \cite{Wright1949}.
Therefore, the stationary density of independent Wright-Fisher diffusions is the product of Dirichlet densities. More precisely,
let
$$\pi(\bar{\x})=\prod_{l=1}^{L}\pi_l(\bar\x^{(l)}), \quad \text{with} \quad
\pi_l(\bar\x^{(l)})=\prod_{i=1}^{M_l-1}(\bar{x}_i^{(l)})^{2 u_i^{(l)}-1}
\left(1-\sum_{i=1}^{M_l-1}\bar{x}_i^{(l)}
\right)^{2 u_{M_l}^{(l)}-1},$$
where $\pi(\bar{\x})$ is the non-normalised stationary density of a coupled Wright-Fisher diffusion with no interaction among loci. In the presence of interaction and assuming parent independent mutations,  i.e. $u_{ij}^{(l)}=u_j^{(l)}$, $\ i,j=1,\dots,M_l, l=1,\dots,L$,
Aurell et al. \cite{timo} prove that there is an additional exponential factor in the stationary density, that is
\begin{equation} \label{eq:stationary}
p(\mathbf{\bar{x}})=\frac{1}{Z}
\pi(\mathbf{\bar{x}}) e^{2 V(\mathbf{\bar{x}})}, 
\end{equation}
with $V$  defined on $\mathcal{\bar{S}}$ by naturally defining the missing frequencies as one minus the sum of the other frequencies at the same locus.
The form of the stationary density is explicit up to a normalising constant. In general, it is difficult to compute the normalising constant $Z$ explicitly,
but under additional assumptions it can be computed numerically, as demonstrated in Section \ref{ASG} and \ref{ex}.

\section{Outline of the derivation of a dual process}
\label{outline}
To derive a process that is dual to the coupled Wright-Fisher diffusion, a generator approach will be used as in \cite{Griffiths2016}, where the authors find a dual process for a multi-locus Wright-Fisher diffusion with recombination. In this section the method will be explained, in general terms.

Let $\mathcal{L}$ be the generator of the diffusion process (\ref{generator}) 
and $\mathcal{L}^D$ be the unknown generator of a dual process.  Suppose that the following relationship holds
\begin{equation}
\label{dualityrelationship}
\mathcal{L}F(\cdot,\n)(\x)=
\mathcal{L}^DF(\x,\cdot)(\n),
\quad
\x\in \mathcal{S},
\quad
\n\in\mathbb{N}^M,
\end{equation}
for some  duality function $F$ that needs to be determined. Using the relationship \eqref{dualityrelationship} the transition rates of a dual process can be identified from its generator. To pursue this approach, it is necessary to compute 
the left hand side of (\ref{dualityrelationship}) by applying the 
generator $\mathcal{L}$  to the duality function $F$,  considered as a function of $\x$,  and rewrite it into the form 
\begin{equation}
\label{LF1}
\mathcal{L}F(\cdot ,\n)(\x)=
\sum_{\hat{\n}}q(\n,\hat{\n})
\left[F(\x,\hat{\n})-F(\x,\n)\right],
\end{equation}
for some non-negative coefficients $q(\n,\hat{\n})$,  $\hat{\n}\in \mathbb{N}^M$, $\hat{\n}\neq \n$. In light of the duality relationship,  expression (\ref{LF1}) can be interpreted as 
the generator $\mathcal{L}^D$  applied to the duality function $F$,  considered as a function of  $\n$. 
Consequently, the dual process obtained in this way is a pure jump Markov process 
on the discrete space $ \mathbb{N}^M $
with transition rate matrix $Q=(q(\cdot,\cdot))$, 
the off-diagonal elements being  the non-negative coefficients in (\ref{LF1}) and the diagonal elements being chosen so that the sum of each row is $0$. 
The alleged duality relationship is validated once  the transition rates and the proper duality function are determined.

Consider the following proposal for the duality  function, $F$. The inspiration for the proposal comes from the duality function for the one-locus Wright-Fisher diffusion with mutations, see e.g. \cite{Etheridge2009,Griffiths2016}. It can be generalised to the multi-locus setting by taking
\begin{equation} \label{eq:dualityfunction}
F(\x,\n)=
\frac{1}{k(\n)}
\prod_{l=1}^{L}\prod_{i=1}^{M_l}
(x_i^{(l)})^{n_i^{(l)}}, 
\end{equation}
for some function $k:\mathbb{N}^M\rightarrow \mathbb{R}\setminus\{0\}$ that is determined in the following. Note that the duality function $F(\cdot,\n)$ defined in \eqref{eq:dualityfunction} belongs to $\mathcal{C}^\infty (\mathcal{S})$, for all $\n\in \mathbb{N}^M$, and thus it belongs to the domain of $\mathcal{L}$. 
Let $\mathbf{\tilde{X}}$ be distributed according to the stationary distribution of the diffusion process $\{\mathbf{X}(t)\}_{t\geq0}$, when such a distribution exists. Then $\E{\mathcal{L}F(\mathbf{\tilde{X}},\n)}=0$.    Therefore, by taking expectation under the stationary distribution in (\ref{LF1}), it follows that
\begin{equation*}
\sum_{\mathbf{\hat{n}}}q(\n,\hat{\n})
\E{F(\mathbf{\tilde{X}},\hat{\n})-F(\mathbf{\tilde{X}},\n)}
=0,
\end{equation*}
which implies that $\E{F(\mathbf{\tilde{X}},\cdot)}$ must be constant. The constant can be taken to be equal to $1$,
and consequently,  
\begin{align}\label{eq:k}
k(\n)=
\mathbb{E}\left[\prod_{l=1}^{L}\prod_{i=1}^{M_l}
(\tilde{X}_i^{(l)})^{n_i^{(l)}}\right].
\end{align}

Note that the existence of a stationary distribution for the diffusion is needed in order to define the function $k$. Thus, in the following,  it is assumed that such a distribution exists. 
Furthermore, in order for the duality function $F$ to be well defined, the function $k$ needs to be non-zero, which holds if
	\begin{equation}
	\label{assumption0}
	\mathbb{P}(\tilde{X}_i^{(l)}=0)=0, \quad i=1,\dots, M_l, l=1,\dots L.
	\end{equation}
In many cases it is possible to verify that a stationary distribution exists and fulfils \eqref{assumption0}.
 For example, as shown in the previous sections,  when the mutations are  parent independent, the stationary density is known, see \eqref{eq:stationary}, and $k(\n)\neq 0 $ for all $\n\in\mathbb{N}^M$.
More generally,
condition \eqref{assumption0} is   satisfied when  the stationary distribution has a density with respect to the Lebesgue measure.
Even if a stationary density is  not known in an explicit form, classical techniques, see e.g. \citep{Khasminskii1980}, may be used to show its existence and properties, using the Fokker-Planck equation in \citep{timo}. 

A relevant case, in which  \eqref{assumption0} is not verified, is the case of no mutations, $\theta=0$. 
Nevertheless, it is still possible to derive a dual process in this case by defining the function $k$ in a simpler  way that does not rely on a stationary distribution. The derivation of the dual process actually becomes simpler than the one outlined in this section. The case of no mutations is treated separately in Section \ref{multi}, Corollary \ref{coroll}. Elsewhere in the paper 
it is assumed that a  stationary distribution  exists and satisfies \eqref{assumption0}.

To  find the transition rates of the dual process, it remains to obtain an expression of the form (\ref{LF1}). 
In fact, it is sufficient to  obtain an expression of the form
\begin{equation}
\label{LF2}
\mathcal{L}F(\cdot , \n)(\x)=
\sum_{\hat{\n}\neq \n}q(\n,\hat{\n})F(\x,\hat{\n})
+q(\n,\n)F(\x,\n),
\end{equation}
with the requirement that $q(\n,\hat{\n})$ is positive for $ \hat{\n}\neq \n$ (it will be soon clear  that $q(\n,\n)$ is thus negative). 
Once  (\ref{LF2}) is obtained, it is possible to derive expression (\ref{LF1}) as follows.
Rewriting  (\ref{LF2}) yields 
\begin{equation}
\label{LF3}
\mathcal{L}F(\cdot , \n)(\x)=
\sum_{\hat{\n}\neq \n}q(\n,\hat{\n})
\left[F(\mathbf{x},\mathbf{\hat{n}})
-F(\mathbf{x},\mathbf{n})\right]
+\Big(\sum_{\hat{\n}\neq \n}q(\n,\hat{\n})+q(\n,\n)\Big)
F(\x,\n).
\end{equation}
Keeping in mind that $\E{\mathcal{L}F(\cdot,\n)(\mathbf{\tilde{X}})}=0$ and that  $\E{F(\mathbf{\tilde{X}},\cdot)}$ is constant, one can apply the expectation with respect to the stationary distribution to get
\begin{equation}
\label{sumq}
\sum_{\hat{\n}\neq \n}q(\n,\hat{\n})+q(\n,\n)=0.
\end{equation}
Therefore (\ref{LF3}) implies (\ref{LF1}) and it remains to write $\mathcal{L}F$ as in (\ref{LF2}) by finding the positive coefficients $q(\n,\hat{\n})$. Furthermore, (\ref{sumq}) can be used to find a recursion formula for the function $k$.

Throughout the rest of the paper, the emphasis will be on obtaining an expression of the type (\ref{LF2}).
This approach is first illustrated in a simpler case (single locus),  in order to lighten the formulas 
and highlight the ideas,  and is subsequently used in the general case of the coupled Wright-Fisher diffusion. 
The simpler case turns out to be closely related to a well known model: the ancestral selection graph. 

\section{The ancestral selection graph}
\label{ASG}
When only one locus is considered, the coupled Wright-Fisher diffusion is simply  a one-locus Wright-Fisher diffusion with selection.
Let $L=1$ , $M_1=2$
and assume that mutations are  parent independent, i.e. $u_{ij}=u_j \; \text{for }  i,j=1,2$. 
The matrix of pairwise selection parameters is the zero matrix and the quadratic function $V$ becomes linear
$$V(\x)=h_1 x_1 +h_2 x_2.$$
Let $j(i)$ be the index opposite to $i$,
$$
j(i)=
\begin{cases}
2 & \text{if } i=1\\
1 & \text{if }  i=2
\end{cases}. 
$$
Then,  the drift terms can be written as follows
\begin{equation*}
\begin{aligned}
& \mu_i(\x)= u_i x _{j(i)} -u_{j(i)}x_i, \\
& g_i(\x) =
h_i x_i(1-x_i) - h_{j(i)} x_i x_{j(i)}, \quad i = 1,2. 
\end{aligned}
\end{equation*}
The diffusion process solving \eqref{SDE} under the assumptions in this section is a two-types Wright-Fisher diffusion with selection and parent independent mutations.
It is known that the genealogical process for this type of Wright-Fisher diffusion is embedded in a graph with coalescing and branching structure, the ancestral selection graph (ASG), studied  by Krone and Neuhauser \cite{Krone1997, Neuhauser1997}. 
In the ASG, first the coalescing-branching structure is constructed  leaving types aside, then types and mutations are superimposed  on it. In contrast, here it is assumed that the types of individuals in the sample $\n$ are known and  mutations are included in the dual process rather than superimposed  afterwards. Our approach is similar to the one in \cite{Etheridge2009}, where the authors  derive a dual process for the finite population size Moran model and use it to find the limiting transition rates of the dual process for the diffusion. 

Following the outline in Section \ref{outline}, a dual process is derived as follows. By applying the generator 
$\mathcal{L}$ to the duality function $F$ in \eqref{eq:dualityfunction}, rewriting the derivatives of $F$, and rearranging the terms yields
\begin{align*}
\mathcal{L}F(\cdot,\mathbf{n})&(\mathbf{x}) =
\\
&
\sum_{i=1,2}
(u_i x _{j(i)} -u_{j(i)}x_i )
\frac{n_i}{x_i}
F(\x,\n)
+
\sum_{i=1,2}
x_i(h_i - h_i x_i - h_{j(i)}x_{ j(i)})
\frac{n_i}{x_i}
F(\x,\n)
\\
& \quad +
\frac{1}{2}
\sum_{i=1,2}
x_i(1-x_i)
\frac{n_i(n_i-1)}{(x_i)^2}
F(\x,\n)
-x_1 x_2 \frac{n_1n_2}{x_1 x_2}F(\x,\n) \\
=&
\sum_{i=1,2}
\frac{n_i(n_i-1)}{2} \frac{1}{x_i} F(\x,\n) 
+
\sum_{i=1,2} 
u_i n_i \frac{x_{j(i)}}{x_i} F(\x,\n) 
-
\sum_{i=1,2}
h_i (n_i +n _{j(i)}) x_i F(\x,\n)\\
&\quad -
\left\{
\frac{n}{2}(n-1) +\sum_{i=1,2}n_i u_{j(i)}-\sum_{i=1,2}n_i h_i
\right\}F(\x,\n), 
\end{align*}
where $n=n_1+n_2$.
To obtain an expression of the form (\ref{LF2}) the expression in the last display can be rewritten as follows.
First replace $x_i=1-x_{j(i)}$ to obtain positive coefficients for the selection terms,  
then use the identities, for $i=1,2$,
\begin{align}
\frac{1}{x_i} F(\x, \n) &= \frac{k(\mathbf{n-e}_i)}{k(\n)} F(\x,\mathbf{n-e}_i),  \label{eq:rewrite1} \\
\frac{x_{j(i)}}{x_i} F(\x,\n) &= \frac{k(\mathbf{n+e}_{j(i)}-\mathbf{e}_i)}{k(\mathbf{n})}
F(\x,\mathbf{n+e}_{j(i)}-\mathbf{e}_i),   \label{eq:rewrite2} \\  x_i F(\x,\n) &= \frac{k(\mathbf{n+e}_i)}{k(\n)}
F(\x,\mathbf{n+e}_i),  \label{eq:rewrite3}
\end{align}
where $\mathbf{e}_i$, $i=1,2$, are the unit vectors  in $\mathbb{N}^2$.
Finally, it yields, 
\begin{equation}
\label{LFasg}
\begin{aligned}
\mathcal{L}F(\cdot,\mathbf{n})(\mathbf{x})=
&\sum_{i=1,2}
\frac{n_i(n_i-1)}{2}
\frac{k(\mathbf{n-e}_i)}{k(\n)} F(\x,\mathbf{n-e}_i) \\
& \quad +
\sum_{i=1,2} 
u_i n_i 
\frac{k(\mathbf{n+e}_{j(i)}-\mathbf{e}_i)}{k(\mathbf{n})}
F(\x,\mathbf{n+e}_{j(i)}-\mathbf{e}_i)
\\
& \quad+
\sum_{i=1,2}
h_{j(i)} n
\frac{k(\mathbf{n+e}_i)}{k(\n)}
F(\x,\mathbf{n+e}_i)
\\
&\quad -
\left\{
\frac{n}{2}(n-1) +\sum_{i=1,2}n_i u_{j(i)} +\sum_{i=1,2} n_{j(i)} h_i
\right\}F(\x,\n),
\end{aligned}
\end{equation}
which is the desired expression. As demonstrated in  Section \ref{outline} the transition rates of a dual process can be identified directly from this expression. Therefore the dual process for the  Wright-Fisher diffusion considered  in this 
section, with respect to  $F$,  is the pure jump Markov process on the state space $\mathbb{N}^2$, with transition rates as follows. The dual process, in state $\n$, jumps to state
\begin{itemize}
	\item  
	$\mathbf{n-e}_i$,  $i=1,2$, s.t. $ n_i\geq 2$, 
	at rate 
	\begin{equation*}
	q(\n,\mathbf{n-e}_i)
	=
	\frac{n_i(n_i-1)}{2}
	\frac{k(\mathbf{n-e}_i)}{k(\n)};
	\end{equation*}
	\begin{flushright}
		\textit{[coalescence]}
	\end{flushright}
	\item 
	$\mathbf{n+e}_{j(i)}-\mathbf{e}_i$,  $i=1,2$, s.t. $ n_i\geq 1$,
	at rate 
	\begin{equation*}
	q(\n,\mathbf{n+e}_{j(i)}-\mathbf{e}_i)
	=
	u_i n_i 
	\frac{k(\mathbf{n+e}_{j(i)}-\mathbf{e}_i)}{k(\mathbf{n})};
	\end{equation*}
	\begin{flushright}
		\textit{[mutation]}
	\end{flushright}
	
	\item  
	$\mathbf{n+e}_i$,  $i=1,2$, 
	at rate 
	\begin{equation*}
	q(\n,\mathbf{n+e}_i)
	=
	h_{j(i)} n
	\frac{k(\mathbf{n+e}_i)}{k(\n)}.
	\end{equation*}
	\begin{flushright}
		\textit{[branching]}
	\end{flushright}
	
\end{itemize}

As anticipated, the dual process just described corresponds to the limiting process in \cite{Etheridge2009}, which is the block counting process of the ancestral selection graph with types and mutations included in the backward evolution.
From the transition rates $q$, it is observed that three types of events are possible for the dual process: mutation, coalescence and branching. The first two  appear also in the Kingman coalescent, while the latter is a virtual addition to the true genealogical process that is characteristic of the ASG. Seen forward in time, branching represents the event that two potential parents are chosen and only the one carrying the advantageous allele reproduces.  Backward in time, when a branching happens, the individual splits into two individuals: its true parent and its virtual (potential) parent. 

To complete the identification of the transition rates, $q(\n,\n)$ is defined as the coefficient of 
$F(\x,\n)$ in (\ref{LFasg}), 
    \begin{equation*}
    q(\n,\n)=-\frac{n}{2}(n-1) -\sum_{i=1,2}n_i u_{j(i)} -\sum_{i=1,2} n_{j(i)} h_i .
    \end{equation*}
The equality (\ref{sumq}) ensures that the sum of each row  of the transition matrix is equal to zero.
Furthermore, by rewriting \eqref{sumq} in terms of the function $k$, a recursion formula can be obtained as in \cite[Thm 5.2]{Krone1997}.   
 The formula, which we omit, is not useful in general to compute $k$ explicitly,  and even in the simpler case of no selection, where the formula, in principle,  could be used, it is  computationally too expensive for practical purposes.
In general it is not possible to find a closed-form expression for $k$ and thus for the transition rates. However, when the mutations are parent independent, as in this example, the stationary density is explicitly known up to a normalising constant $Z$ and thus  $k$ can be written as an integral with respect to the stationary density
    \begin{equation*}
    k(\n)= \frac{1}{Z} 
    \int_0^1  \!\!\!\!
    x^{n_1+ 2u_1-1}(1-x)^{n_2+2u_2-1} 
    e^{ 2[h_1 x  + h_2  (1-x)]}
    dx.
    \end{equation*}
The integral above cannot be computed analytically but it is related to the confluent hypergeometric function of the first kind, the Kummer function, which can be efficiently computed numerically. 
The idea of using the Kummer function originates from \cite{timo} and \citep{Krone1997}, where it has been used to find, respectively, a series representation for the normalising constant and a representation for the expected allele frequency.  Let ${}_{1}F_1$ be the confluent hypergeometric function, then, using its integral representation, it yields 
\begin{align*}
k(\n) &= 
\frac{1}{Z} e^{2 h_2}
B(n_1+2u_1,n_2+2u_2)
{}_{1}F_1\left(n_1+2u_1,n_1+2u_1+n_2+2u_2,2(h_1-h_2)\right),  \\
Z& = e^{2 h_2}
B(2u_1,2u_2)
{}_{1}F_1(2u_1 , 2u_1 + 2u_2 , 2(h_1-h_2)),
\end{align*}
where $B$ is the Beta function. See \cite{Abramowitz} for a complete collection of definitions and properties of confluent hypergeometric functions.

\section{A multi-locus dual process}
\label{multi}
In this section  a dual process for the coupled Wright-Fisher diffusion is derived in the general  multi-locus setting, $L \geq 1$ and $M_l \geq 2$, $l=1, \dots, L$.  

\begin{theorem}\label{thm:main}
Let $\X$ be the coupled Wright-Fisher diffusion  with generator \eqref{generator}, where $\mu$, $g$ and $d$ are given by \eqref{mu}, \eqref{g} and \eqref{d}, respectively.
Assume  a stationary distribution for the diffusion exists and satisfies \eqref{assumption0}.
Let $k$ be given by \eqref{eq:k} and let the duality function $F$ be given by \eqref{eq:dualityfunction}.
Then there exists a dual process $\mathbf{N}$ for $\X$, in the sense of  \eqref{duality_expectations}, with respect to the duality function $F$, where $\mathbf{N}$  is the pure jump Markov process on the state space $\mathbb{N}^M$ with the following transition rates. From the current state, $\n 
\in \mathbb{N}^M\setminus\{0\}$, $\mathbf{N}$ jumps to
	\begin{itemize}
		\item $\n - \e_i^{(l)}, i=1,\dots M_l, l=1,\dots, L,$ s.t. $n_i^{(l)}\geq 2$,
		at rate
		$$
		q(\n, \n - \e_i^{(l)} )=
		\frac{n_i^{(l)}(n_i^{(l)}-1)}{2}
		\frac{k(\n - \e_i^{(l)})}{k(\n)} ;
		$$
		\begin{flushright}
		[coalescence]
		\end{flushright}
		\newpage
		\item $\n - \e_i^{(l)} + \e_j^{(l)} , i,j=1,\dots M_l, l=1,\dots, L,$ s.t. $i\neq j,n_i^{(l)}\geq 1$,
		at rate
		$$
		q(\n, \n - \e_i^{(l)}+ \e_j^{(l)}  )=
		n_i^{(l)}u_{ji}^{(l)}
		\frac{k(\n - \e_i^{(l)}+ \e_j^{(l)} )}{k(\n)} ;
		$$
		\begin{flushright}
			[mutation]
		\end{flushright}
		
		\item $\n + \e_j^{(l)}, j=1,\dots M_l, l=1,\dots, L,$
		at rate
		$$
		q(\n, \n + \e_j^{(l)} )=
		\left( 
		n^{(l)}\sum_{\substack{k=1\\k\neq j}}^{M_l}h_{k}^{(l)} 
		+\sum_{\substack{r=1\\r\neq l}}^{L}\sum_{i=1}^{M_r}n_i^{(r)}J_{ji}^{(lr)}
		\right)
		\frac{k(\n + \e_j^{(l)})}{k(\n)} ;
		$$
		\begin{flushright}
			[single-branching]
		\end{flushright}
		
		\item $\n + \e_j^{(l)}+\e_h^{(r)}, j=1\dots M_l, h=1\dots,M_r,  l,r=1,\dots, L, r>l,$ 
		at rate
		$$
		q(\n, \n + \e_j^{(l)}+\e_h^{(r)} )=
		\left( (n^{(l)}+n^{(r)})
		\sum_{\substack{k=1\\k\neq j}}^{M_l}\sum_{\substack{m=1\\m\neq h}}^{M_r}	
		J_{km}^{(lr)}	\right)
		\frac{k(\n + \e_j^{(l)}+\e_h^{(r)})}{k(\n)} ; 
		$$
		\begin{flushright}
			[double-branching]	
		\end{flushright}
		where $n^{(l)}=\norm{\n^{(l)}}=\sum_{i=1}^{M_l}n_i^{(l)}$. Furthermore,
		\begin{equation*}
		q(\n,\n)
		= - \sum_{l=1}^{L}  
		\left(
		\frac{1}{2}n^{(l)}(n^{(l)}-1) 
		+
		\sum_{\substack{i,j=1\\i\neq j}}^{M_l}
	n_i^{(l)}u_{ij}^{(l)} +
		\sum_{i=1}^{M_l}\sum_{\substack{k=1\\k\neq i}}^{M_l}
   	 	n_{i}^{(l)}h_{k}^{(l)}
    		+\sum_{\substack{r=1\\r\neq l}}^{L}\sum_{k=1}^{M_l}\sum_{m=1}^{M_r}	
    		n^{(l)}J_{km}^{(lr)} 
		\right).
		\end{equation*}
	\end{itemize}
\end{theorem}

Note that the mutation and coalescence jumps involve one locus at the time. The coalescence and mutation  rates are  similar  to the transition rates of the  Kingman coalescent process with mutations, the only difference being the function $k$ which, despite having the same structure, is based on a different stationary density and depends on all the loci, not only on the one where the jump takes place.
The single-branching rate  not only contains the single-locus selection parameters  in a form that generalises the rates in Section \ref{ASG}, but it also contains the two-locus selection parameters to include the effect of the pairwise interaction on the single locus. The single-branching also involve only one locus at the time. Finally, the double-branching rate reflects the particular structure of pairwise interactions of the coupled Wright-Fisher diffusion and it is, to the best of our knowledge,  a novel type of transition rate appearing in  genealogical processes related to Wright-Fisher diffusions. 
The double-branching represents simultaneous branching at two different loci. 
As anticipated in the introduction, the dual process can  thus be interpreted as the block counting process of $L$ coupled ASGs.
Furthermore, if $J=0$,  the loci are independent since $\nabla V=\mathbf{h}$, and thus double-branching does not occur and the dual
process consists of $L$ independent ASGs, as in \cite{Etheridge2009} with $-2 \sum_{\substack{k=1\\k\neq j}}^{M_l}h_{k}^{(l)} $ corresponding to $\sigma_j $.

The explicit parts of the transition rates, not depending on the function $k$, have a very natural interpretation. As in the simpler case studied in Section \ref{ASG}, the basic principle is that weak  types branch at a higher rate. The difference is that, while in the simpler case there are only two types, a viable type and a weaker type, here  there are many types and many loci all influencing each other's branching rates. 
To understand this behaviour in greater detail,  some terms will be investigated more thoroughly.
The term
$$
n^{(l)}\sum_{\substack{k=1\\k\neq j}}^{M_l}h_{k}^{(l)},  
$$
arises purely from the one-locus selection  and contributes to the rate of adding a gene of type $j$ at locus $l$. It depends on the one-locus viability of the other allele types  (all except type $j$) at locus $l$, the higher their viability, the higher the rate of adding  type $j$, and of course it is also directly proportional to the number, $n^{(l)}$, of genes at locus $l$.

The rate of adding a couple of genes of type $j$ at locus $l$ and of type $h$ at locus $r$ is related  to 
the term 
    $$ 
    (n^{(l)}+n^{(r)})
    \sum_{\substack{k=1\\k\neq j}}^{M_l}\sum_{\substack{m=1\\m\neq h}}^{M_r}
    J_{km}^{(lr)}.
    $$
It depends on the  viability of the other  couples of allele types (all except couple $j,h$) at loci $l$ and $r$, the higher their viability, the higher the rate of adding type $j$ and $h$ at locus $l$ and  $r$,  respectively.
Again the rate is  directly proportional to the number, $n^{(l)}+n^{(r)}$,  of genes at loci $l$ and $r$.

Although the interpretation of some parts of the transition rates is straightforward, the function $k$ remains implicit, similar to the simpler Kingman coalescent process and the ancestral selection graph with parent dependent mutations. When the mutations are parent independent, the stationary density is known up to a normalising constant and $k$ can be expressed as an integral that sometimes can be easily  computed numerically, see Section \ref{ASG} and Section \ref{ex}, where a series representation of $k$ involving Kummer and Beta functions will be given.
Nevertheless, even when the stationary distribution is not explicitly known, but still exists,  Theorem \ref{thm:main} provides information on the structure of the transition rates of the dual process that may be useful.  As explained in the introduction, many established inference methods for populations under various generalisations of the Wright-Fisher diffusion rely on approximating the backward dynamics of the associated genealogical process. Deriving a dual process for the coupled Wright-Fisher diffusion is central to further investigations concerning the genealogy of a sample and possibly provides a basis for the construction of inference methods inspired by the existing methods described in the introduction.

In general, if the transition rates are not known explicitly, it might seem difficult to provide bounds for the dual process. However, its growth is controlled by a much simpler Markov chain. Indeed, the process  $ \norm{\N(t)}$ is a jump process on  $\mathbb{N}\setminus\{0\}$, with possible jumps: $+2$, representing double-branching, $+1$, representing  single-branching and $-1$,  representing coalescence.
As is typical of genealogical processes appearing in population genetics, 
 the rate of negative (coalescent) jumps is at least quadratic and the rate of positive (branching) jumps is at most linear, as shown in details in \eqref{ratesbounds} in the Appendix. This allows to construct a monotone coupling to bound  the jump chain of $\norm{\N}$ by a simpler Markov chain  which is reported  in the Appendix as it could be useful for future work.

This section concludes with the extension of Theorem \ref{thm:main} to the case of no mutations, Corollary \ref{coroll}, and with two applications of the duality relationship. The first one is useful to derive an expansion of the transition density of the diffusion and the second one to study, in the absence of mutations, fixation/extinction probabilities of allele types. 

The duality relation (\ref{duality_expectations}), which follows from  Theorem \ref{thm:main}, can be rewritten as 
$$
\E{S(\X(t),\n)|\X(0)=\x}= \E{S(\tilde{\X},\n)}\sum_{\m\in \mathbb{N}^M} p_{\n,\m}(t)F(\x,\m),
$$
where $S(\x,\n)$ is the probability mass function of $L$ independent multinomial random variables with parameters $\x^{(1)},\dots,\x^{(L)}$ and $n^{(1)},\dots,n^{(L)}$, and $p_{\n,\m}(t)$ are the transition densities of the dual process $\N$.
By applying  \textit{sample inversion}  as $n\to\infty$, as in \cite{Etheridge2009,Griffiths2016}, an expansion for the transition density of $\X$ can be obtained in terms of the limit of the transition densities of $\N$, the stationary density  of $\X$ and the duality function $F$. This corresponds to identifying the distribution of $\X(t)$ from its moments. The derivation of the expansion is essentially similar to the one in \cite{Griffiths2016}, a rigorous proof is left to future work.

In the case of no mutation, $\theta_l=0$, $l=1,\dots,L$, the boundaries are absorbing for the diffusion $\X$.  
Any distribution that puts all its mass on one allele type for each locus is an invariant distribution for the diffusion but does not satisfy 
assumption \eqref{assumption0}. Nevertheless, 
as anticipated in Section \ref{outline}, 
 the derivation of  a dual process in this case is simpler than in the presence of mutations,
as there is no need of relying on  invariant distributions to define the duality function.
In fact, the duality function can be defined explicitly by defining  the function $k$ to be equal to a product of multivariate Beta functions,
\begin{equation}
\label{k_nomut}
k(\n)=\prod_{l=1}^L B(\n^{(l)})
\quad \text{with} \quad
B(\n^{(l)})=\frac{\prod_{i=1}^{M_l} \Gamma(n_i^{(l)})}{\Gamma(n^{(l)})}
.
\end{equation}
To get an intuition on why $k$  is defined in this way, note that 
in the neutral one-locus model with parent independent mutations 
$k(\n)=\frac{B(\n+2\mathbf{u})}{B(2  \mathbf{u})}$, where $\mathbf{u}$ is the vector of mutation transition rates with $u_j=\theta P_{ij}$,
and, as $\theta\to 0$, $k(\n)$ converges to  $B(\n)$ .
Using definition \eqref{k_nomut} for $k$, the transition rates of the dual process can be derived from those in Theorem \ref{thm:main}, see Section \ref{proof} for more details, to obtain the following 
\begin{corollary}
\label{coroll}
Let $\X$ be defined as in Theorem \ref{thm:main} and assume $\theta_l=0$, $l=1,\dots,L$. 
Then there exists a dual process $\mathbf{N}$ for $\X$, in the sense of  \eqref{duality_expectations}, with respect to the duality function $F(\x,\n)=\prod_{l=1}^L\frac{1}{B(\n^{(l)})}\prod_{i=1}^{M_l} (x_{i}^{(l)})^{n_i^{(l)}}$.
 $\mathbf{N}$  is the pure jump Markov process on the state space $\mathbb{N}^M$ with  transition rates
\begin{itemize}
\item 	$
		q(\n, \n - \e_i^{(l)} )=
		\frac{n_i^{(l)}(n^{(l)}-1)}{2} ,
		  \quad i=1,\dots M_l, l=1,\dots, L ;
		$
\item 
		$
		q(\n, \n + \e_j^{(l)} )=
		n_j^{(l)}\sum_{\substack{k=1\\k\neq j}}^{M_l}h_{k}^{(l)} 
		+\frac{n_j^{(l)}}{n^{(l)}}\sum_{\substack{r=1\\r\neq l}}^{L}\sum_{i=1}^{M_r}n_i^{(r)}J_{ji}^{(lr)} , \quad
		j=1,\dots M_l, l=1,\dots, L;
		$
\item	$
		q(\n, \n + \e_j^{(l)}+\e_h^{(r)} )=
		 (n^{(l)}+n^{(r)}) \frac{n_j^{(l)}n_h^{(r)}}{n^{(l)}n^{(r)}}
		\sum_{\substack{k=1\\k\neq j}}^{M_l}\sum_{\substack{m=1\\m\neq h}}^{M_r}	
		J_{km}^{(lr)}
		, \quad  j=1\dots M_l, h=1\dots,M_r,$ $  l,r=1,\dots, L, r>l.
		$
\end{itemize}
\end{corollary}
Note that, due to the absence of mutations, the transition rates in this case are fully explicit. 
By Corollary \ref{coroll}, the duality relationship \eqref{duality_expectations} can be rewritten  as 
\begin{equation}
\E{\prod_{l=1}^L\prod_{i=1}^{M_l} X_{i}^{(l)}(t)^{n_i^{(l)}}|\X(0)=\x}=
\E{\prod_{l=1}^L \frac{B(\n^{(l)})}{B(\N^{(l)}(t))} \prod_{i=1}^{M_l} (x_{i}^{(l)})^{N_i^{(l)}(t)} |\N(0)=\n} ,    
\end{equation}
letting $t \to \infty$ enables the study of fixation/extinction probabilities of allele types, as  in e.g. \cite{Etheridge2009,Griffiths2016,Mano2009,Foucart2013,gonzalez2018}.

\section{Proofs of the main results }
\label{proof}
\subsection{Proof of Theorem \ref{thm:main}}
Following the outline in Section \ref{outline}, a dual process is derived as follows. 
By applying the generator 
$\mathcal{L}$ to the duality function $F$ in \eqref{eq:dualityfunction} each term in the expression for $\mathcal{L}F$ is treated separately. As in  Section \ref{ASG}, the terms arising from mutation and diffusion can be easily rewritten in the required form. Summing the mutation terms over allele types at locus $l$ yields
\begin{equation*}
\begin{aligned}
\sum_{i=1}^{M_l}
\mu_i^{(l)}(\x^{(l)})
\frac{\partial F}{\partial x_i^{(l)} }(\mathbf{x},\mathbf{n})
=
\sum_{\substack{i,j=1\\i\neq j}}^{M_l}
n_i^{(l)}u_{ji}^{(l)} \frac{x_j^{(l)}}{x_i^{(l)}}F(\mathbf{x},\mathbf{n})
-
\sum_{\substack{i,j=1\\i\neq j}}^{M_l}
n_i^{(l)}u_{ij}^{(l)}
F(\mathbf{x},\mathbf{n}).
\end{aligned}
\end{equation*}
Using identity (\ref{eq:rewrite2}) at locus $l$ the mutation terms can be rewritten in the desired form 
\begin{equation}
\begin{aligned}
\label{mutation}
\sum_{i=1}^{M_l}
\mu_i^{(l)}(\x^{(l)})
&\frac{\partial F}{\partial x_i^{(l)} }(\mathbf{x},\mathbf{n})  \\
=
&\sum_{\substack{i,j=1\\i\neq j}}^{M_l}
n_i^{(l)}u_{ji}^{(l)} 
\frac{k(\n-\e_i^{(l)}+\e_j^{(l)})}{k(\n)} F(\x,\n-e_i^{(l)}+e_j^{(l)})
-
\sum_{\substack{i,j=1\\i\neq j}}^{M_l}
n_i^{(l)}u_{ij}^{(l)}
F(\mathbf{x},\mathbf{n}).
\end{aligned}
\end{equation}
For the diffusion terms, the diagonal and off-diagonal terms are written separately as
\begin{align*}
d_{ii}^{(l)} (\x^{(l)})  \frac{\partial^2F}{\partial x_i^{(l) 2}  } (\x,\n) &=	
n_i^{(l)}(n_i^{(l)}-1) \frac{1}{x_i^{(l)}} F(\x,\n) - n_i^{(l)}(n_i^{(l)}-1)F(\x,\n),
\\ 
d_{ij}^{(l)} (\x^{(l)})\frac{\partial^2F}{\partial x_i^{(l)} \partial x_j^{(l)} }(\x,\n) & =
-n_i^{(l)}n_j^{(l)}F(\x,\n),
\quad\quad i \neq j.
\end{align*}
Summing  the diffusion terms at locus $l$ and rearranging yields
\begin{equation*}
\begin{aligned} 
\frac{1}{2}\sum_{i,j=1}^{M_l}
d_{ij}^{(l)}(\x^{(l)})
\frac{\partial^2F}{\partial x_i^{(l)} \partial x_j^{(l)} }
(\mathbf{x},\mathbf{n})=&
\sum_{i=1}^{M_l}
\frac{n_i^{(l)}(n_i^{(l)}-1)}{2} \frac{1}{x_i^{(l)}} F(\x,\n)
-
\frac{1}{2}n^{(l)}(n^{(l)}-1)F(\x,\n).
\end{aligned}
\end{equation*}	
Now use identity (\ref{eq:rewrite1}) at locus $l$ to obtain 
\begin{equation}
\label{diffusion}
\begin{aligned}
\frac{1}{2}\sum_{i,j=1}^{M_l}
d_{ij}^{(l)}(\x^{(l)})
\frac{\partial^2F}{\partial x_i^{(l)} \partial x_j^{(l)} }
(\mathbf{x},\mathbf{n})= &
\sum_{i=1}^{M_l}
\frac{n_i^{(l)}(n_i^{(l)}-1)}{2} 
\frac{k(\n-\e_i^{(l)})}{k(\n)} F(\x,\n-\e_i^{(l)})\\
& -
\frac{1}{2}n^{(l)}(n^{(l)}-1)F(\x,\n).
\end{aligned}
\end{equation}	
Next, consider  the interaction terms.  Using the definition of $g^{(l)}$ and rewriting the derivatives of $F$ yields, 
\begin{equation}
\label{sumgil}
\begin{aligned}
\sum_{i=1}^{M_l}
g_i^{(l)}(\mathbf{x})
\frac{\partial F}{\partial x_i^{(l)} }(\mathbf{x},\mathbf{n})
= &
\underbrace{
	\sum_{i=1}^{M_l}
	\sum_{k=1}^{M_l}
	n_{i}^{(l)}h_{k}^{(l)}(\delta_{ik}-x_{k}^{(l)})F(\x,\n)
}_{S_1}\\
&\quad +
\underbrace{
	\sum_{\substack{r=1\\r\neq l}}^{L}\sum_{i=1}^{M_l}\sum_{k=1}^{M_l}\sum_{m=1}^{M_r}	
	n_{i}^{(l)}J_{km}^{(lr)}(\delta_{ik}-x_{k}^{(l)})x_m^{(r)}F(\x,\n)
}_{S_2}.
\end{aligned}
\end{equation}
Note that the first group of sums, $S_1$, contains the one-locus selection parameters while the second, $S_2$, contains the pairwise selection parameters. Each of them will be treated separately.
The one-locus selection term can be rearranged into
\begin{equation*}
\begin{aligned}
S_1=
\sum_{i=1}^{M_l}
n_{i}^{(l)}h_{i}^{(l)}F(\x,\n)
-
\sum_{i=1}^{M_l}
\sum_{k=1}^{M_l}
n_{i}^{(l)}h_{k}^{(l)}x_{k}^{(l)}F(\x,\n).
\end{aligned}
\end{equation*}
As in Section \ref{ASG}, the fact that the sum of the frequencies at each locus equals one is used. Since  $x_{k}^{(l)}=1- \sum_{\substack{j=1\\j\neq k}}^{M_l}x_{j}^{(l)}$, the terms can be rearranged to obtain
\begin{equation*}
\begin{aligned}
S_1=
-\sum_{i=1}^{M_l}
\sum_{\substack{k=1\\k\neq i}}^{M_l}
n_{i}^{(l)}h_{k}^{(l)}F(\x,\n)
+
\sum_{j=1}^{M_l}
\left( n^{(l)}\sum_{\substack{k=1\\k\neq j}}^{M_l}h_{k}^{(l)} \right)
x_{j}^{(l)}F(\x,\n).
\end{aligned}
\end{equation*}
The second part of \eqref{sumgil} can be expressed as
\begin{equation*}
\begin{aligned}
S_2
=
\sum_{\substack{r=1\\r\neq l}}^{L}\sum_{i=1}^{M_l}\sum_{m=1}^{M_r}	
n_{i}^{(l)}J_{im}^{(lr)}x_m^{(r)}F(\x,\n)
-
\sum_{\substack{r=1\\r\neq l}}^{L}\sum_{k=1}^{M_l}\sum_{m=1}^{M_r}	
n^{(l)}J_{km}^{(lr)}x_{k}^{(l)}x_m^{(r)}F(\x,\n).
\end{aligned}
\end{equation*}
This time the equality 
\begin{equation}
\label{trick}
-x_k^{(l)}x_m^{(r)}=
-1+ 
\sum_{j=1}^{M_l}
\sum_{h=1}^{M_r}
(1-\delta_{hm}\delta_{jk})x_j^{(l)}x_h^{(r)},
\end{equation}
will be used.  
To see that (\ref{trick}) holds, the fact that the frequencies sum up to one at each locus is used multiple times, as 
follows, 
\begin{equation*}
\begin{aligned}
-x_k^{(l)}x_m^{(r)}
=&
-x_k^{(l)}\left(
1-\sum_{\substack{h=1\\h\neq m}}^{M_r}x_h^{(r)}
\right)\\
=&
-x_k^{(l)} + 
\sum_{h=1}^{M_r}(1-\delta_{hm})x_k^{(l)}x_h^{(r)}\\
=&
-1 + \sum_{\substack{j=1\\j\neq k}}^{M_l}x_j^{(l)}
\cdot\sum_{h=1}^{M_r}x_h^{(r)}
+\sum_{h=1}^{M_r}(1-\delta_{hm})x_k^{(l)}x_h^{(r)}\\
=&
-1 
+\sum_{j=1}^{M_l}\sum_{h=1}^{M_r}
(1-\delta_{jk})	x_j^{(l)}x_h^{(r)}
+\sum_{j=1}^{M_l}\sum_{h=1}^{M_r}
\delta_{jk}(1-\delta_{hm})	x_k^{(l)}x_h^{(r)}\\
=&
-1+ 
\sum_{j=1}^{M_l}
\sum_{h=1}^{M_r}
(1-\delta_{hm}\delta_{jk})x_j^{(l)}x_h^{(r)}.
\end{aligned}
\end{equation*}
Applying (\ref{trick}) in the expression for $S_2$ and rearranging the terms, leads to
\begin{equation*}
\begin{aligned}
S_2
=&
-
\sum_{\substack{r=1\\r\neq l}}^{L}\sum_{k=1}^{M_l}\sum_{m=1}^{M_r}	
n^{(l)}J_{km}^{(lr)}F(\x,\n)
+
\sum_{\substack{r=1\\r\neq l}}^{L}\sum_{m=1}^{M_r}	
\left(\sum_{i=1}^{M_l} n_{i}^{(l)}J_{im}^{(lr)}\right)
x_m^{(r)}F(\x,\n)\\
&\quad +
\sum_{\substack{r=1\\r\neq l}}^{L}\sum_{j=1}^{M_l}\sum_{h=1}^{M_r}	
\left(
n^{(l)}\sum_{\substack{k=1\\k\neq j}}^{M_l}\sum_{\substack{m=1\\m\neq h}}^{M_r}	
J_{km}^{(lr)}
\right)
x_{j}^{(l)}x_h^{(r)}F(\x,\n).
\end{aligned}
\end{equation*}
Summing over $l$ and putting similar terms together yields
\begin{equation*}
\begin{aligned}
\sum_{l=1}^{L}
\sum_{i=1}^{M_l}
g_i^{(l)}(\mathbf{x})
\frac{\partial F}{\partial x_i^{(l)} }(\x,\n)
& =
-
\left(
\sum_{l=1}^{L}\sum_{i=1}^{M_l}\sum_{\substack{k=1\\k\neq i}}^{M_l}
n_{i}^{(l)}h_{k}^{(l)}
+\sum_{l=1}^{L}\sum_{\substack{r=1\\r\neq l}}^{L}\sum_{k=1}^{M_l}\sum_{m=1}^{M_r}	
n^{(l)}J_{km}^{(lr)}
\right)F(\x,\n)	
\\
& \quad +	
\sum_{l=1}^{L}\sum_{j=1}^{M_l}
\left( 
n^{(l)}\sum_{\substack{k=1\\k\neq j}}^{M_l}h_{k}^{(l)} 
+\sum_{\substack{r=1\\r\neq l}}^{L}\sum_{i=1}^{M_r}n_i^{(r)}J_{ji}^{(lr)}
\right)
x_{j}^{(l)}F(\x,\n)
\\ 
&\quad +
\sum_{l=1}^{L}\sum_{\substack{r=1\\r> l}}^{L}\sum_{j=1}^{M_l}\sum_{h=1}^{M_r}	
\left(
(n^{(l)}+n^{(r)})
\sum_{\substack{k=1\\k\neq j}}^{M_l}\sum_{\substack{m=1\\m\neq h}}^{M_r}	
J_{km}^{(lr)}
\right)
x_{j}^{(l)}x_h^{(r)}F(\x,\n).
\end{aligned}
\end{equation*}
Use the  identities (\ref{eq:rewrite3}) at locus $l$ and
    \begin{equation*}
    \begin{aligned}
    x_{j}^{(l)}x_h^{(r)}F(\x,\n)=
    \frac{k(\n+\e_j^{(l)}+\e_h^{(r)})}{k(\n)} F(\x,\n+\e_j^{(l)}+\e_h^{(r)})
    \end{aligned}
    \end{equation*}
for the mixed terms involving loci $l$ and $r$, in order to rewrite the selection terms in the desired form
    \begin{equation}
    \label{selection}
    \begin{aligned}
    \sum_{l=1}^{L}
    &
    \sum_{i=1}^{M_l}
    g_i^{(l)}(\mathbf{x})
    \frac{\partial F}{\partial x_i^{(l)} }(\x,\n)
    =\\
    & -
    \left(
    \sum_{l=1}^{L}\sum_{i=1}^{M_l}\sum_{\substack{k=1\\k\neq i}}^{M_l}
    n_{i}^{(l)}h_{k}^{(l)}
    +\sum_{l=1}^{L}\sum_{\substack{r=1\\r\neq l}}^{L}\sum_{k=1}^{M_l}\sum_{m=1}^{M_r}	
    n^{(l)}J_{km}^{(lr)}
    \right)F(\x,\n)	
    \\
    & +	
    \sum_{l=1}^{L}\sum_{j=1}^{M_l}
    \left( 
    n^{(l)}\sum_{\substack{k=1\\k\neq j}}^{M_l}h_{k}^{(l)} 
    +\sum_{\substack{r=1\\r\neq l}}^{L}\sum_{i=1}^{M_r}n_i^{(r)}J_{ji}^{(lr)}
    \right)
    \frac{k(\n+\e_j^{(l)})}{k(\n)} F(\x,\n+\e_j^{(l)})
    \\ 
    & +
    \sum_{l=1}^{L}\sum_{\substack{r=1\\r> l}}^{L}\sum_{j=1}^{M_l}\sum_{h=1}^{M_r}
    \left(
    (n^{(l)}+n^{(r)})
    \sum_{\substack{k=1\\k\neq j}}^{M_l}\sum_{\substack{m=1\\m\neq h}}^{M_r}	
    J_{km}^{(lr)}
    \right)
    \frac{k(\n+\e_j^{(l)}+\e_h^{(r)})}{k(\n)} F(\x,\n+\e_j^{(l)}+\e_h^{(r)}).
    \end{aligned}
    \end{equation}
The terms arising from mutation (\ref{mutation}), diffusion (\ref{diffusion}) and selection (\ref{selection}) are now  written in form \eqref{LF2}. It is finally possible  to identify the transition rates of the dual process.

In order to complete the proof, the method of duality is applied, more precisely, Corollary 4.4.13 in \cite{Ethier1986}, which amounts to verifying the following conditions: 
\begin{align*}
     F(\X(t),\n) -\int_0^t \mathcal{L}F(\cdot,\n)(\X(s))ds \quad \text{ and } \quad 
 F(\x,\N(t)) - \int_0^t \mathcal{L}^D F(\x,\cdot)(\N(s))ds
\end{align*}
are martingales and, for each $T > 0$, there exists an integrable random variable $\Gamma_T$ such that
	\begin{align}
	\label{intcond1}
	&\sup_{0\leq s,t\leq T} \left|  F(\X(s),\N(t))\right| \leq \Gamma_T, \\
	&\sup_{0\leq s,t\leq T} \left|  \mathcal{L}F(\cdot,\N(t))(\X(s))\right| 
	=\sup_{0\leq s,t\leq T} \left|  \mathcal{L}^D F(\X(s),\cdot)(\N(t))\right| 
	\leq \Gamma_T,\label{intcond2}  
	\end{align}
almost surely.
First note that,  
 as discussed in Section \ref{outline}, $F(\cdot,\n)$ belongs to the domain of $\mathcal{L}$, for all $\n\in \mathbb{N}^M$, and the integrability conditions 
 \eqref{intcond1}-\eqref{intcond2} ensure 
 that  $F(\x,\cdot)$ belongs to the domain of  $\mathcal{L}^D$, for all $\x\in \mathcal{S}$, 
 and that the processes 
 $F(\X(t),\n) -\int_0^t \mathcal{L}F(\cdot,\n)(\X(s))ds$ and 
$F(\x,\N(t)) -\int_0^t \mathcal{L}^D F(\x,\cdot)(\N(s))ds$ 
are integrable. The martingale property  trivially  holds for both processes.

In order to complete the proof,  
following \cite{barbour2000}  to verify \eqref{intcond1}-\eqref{intcond2}, it is sufficient to show that  there exists a function $H: \mathbb{N}^M\to [0,\infty)$ such that
	\begin{equation}
	\label{fctH}
	F(\x,\n) + \left| \mathcal{L^D}F(\x,\cdot)(\n)\right| \leq H(\n), \quad \forall (\x,\n)\in \mathcal{S}\times 	\mathbb{N}^M,
 	\end{equation}
and $\left\{  H(\N(t\wedge \tau_j)), 0\leq t \leq T,  j\geq 1\right\} $, where 
$\tau_j=\inf \{s\geq 0 : \norm{\N(s)}\geq j \}$,
 is uniformly integrable, for all initial conditions, $\N(0)=\n \in \mathbb{N}^M$, and all $T\geq 0$.
First, bounds for $F$ and $\mathcal{L}^D F$ are provided.

The definition \eqref{eq:k} of $k$ and Jensen's inequality yield 
    $$
	k(\n)=
	\mathbb{E}\left[\prod_{l=1}^{L}\prod_{i=1}^{M_l}
    (\tilde{X}_i^{(l)})^{n_i^{(l)}}\right]
    \geq 
    \mathbb{E}\left[\left(\prod_{l=1}^{L}\prod_{i=1}^{M_l}
    \tilde{X}_i^{(l)}\right)^{\norm{\n}}\right]
    \geq
    \mathbb{E}\left[\prod_{l=1}^{L}\prod_{i=1}^{M_l}
    \tilde{X}_i^{(l)}\right]^{\norm{\n}}.
	$$
Denote $\mathbb{E}\left[\prod_{l=1}^{L}\prod_{i=1}^{M_l}
    \tilde{X}_i^{(l)}\right]^{-1}$ by $a$. Because of assumption \eqref{assumption0}, the latter expectation  is non-zero and  $a$ is well defined,  furthermore, $a>1$.
Consequently, using the definition \eqref{eq:dualityfunction} of $F$,  it follows that
	\begin{equation}
	\label{boundF}
	\begin{aligned}
	F(\x,\n)	\leq a ^{\norm{\n}}.
	\end{aligned}
	\end{equation}
Moreover, 
	\begin{align}
	\left|\mathcal{L} F(\cdot,\n)(\x) \right| =
	\left|\mathcal{L}^D F(\x,\cdot)(\n) \right| 
	&\leq
	\sum_{\hat{\n}} |q(\n,\hat{\n})| |F(\x,\hat{\n})-F(\x,\n)| \nonumber 
	\\
	&\leq
	\sum_{\hat{\n}} |q(\n,\hat{\n})| 
	\left[
	a ^{\norm{\n}+2}	+ a ^{\norm{\n}}
	\right] \nonumber 
	\\
	&\leq
	b'' \norm{\n}^2 a ^{\norm{\n}} ,	\label{boundLF}
	\end{align}
since 
$\norm{\hat{\n}}\leq \norm{\n}+2$, 
$\sum_{\hat{\n}} |q(\n,\hat{\n})|= 2 |q(\n,\n)|$, 
and 
$|q(\n,\n)|\leq b''' \norm{\n}^2$, for some $b'''>0$ and $b''=4b'''a^2$.
By \eqref{boundF} and \eqref{boundLF}, it is implied that  inequality \eqref{fctH} holds true with $H(\n)= b' \norm{\n}^2 a ^{\norm{\n}}$, for some $b'>0$.
Consider now the following representation formula for the expectation of $H$ applied to the stopped dual process,
	\begin{equation}
	\label{UI}
	\E{H(\N(t\wedge\tau_j))| \N(0)=\n}=
	H(\n) +\E{\int_0^{t\wedge \tau_j} 
	\mathcal{L}^D H (\N(s)) ds}.
	\end{equation}
Using the 
inequalities \eqref{ratesbounds} of  the Appendix,
	\begin{equation*}
	\begin{aligned}
	\mathcal{L}^D H (\n)&=  b'
	\sum_{l=1}^{L}\sum_{i=1}^{M_l} q(\n,\n-\e_i^{(l)})
	\left[
	(\norm{\n}-1)^2 a ^{\norm{\n}-1} - \norm{\n}^2 a ^{\norm{\n}}
	\right]
	\\ & + b'
	\sum_{l=1}^{L}\sum_{i=1}^{M_l} q(\n,\n+\e_i^{(l)})
	\left[
	(\norm{\n}+1)^2 a ^{\norm{\n}+1} - \norm{\n}^2 a ^{\norm{\n}}
	\right]
	\\ & + b'
	\sum_{l=1}^{L}\sum_{\substack{r=1\\r> l}}^{L}\sum_{j=1}^{M_l}\sum_{h=1}^{M_r}
	q(\n,\n+\e_j^{(l)}+\e_h^{(r)})
	\left[
	(\norm{\n}+2)^2 a ^{\norm{\n}+2} -\norm{\n}^2 a ^{\norm{\n}}
	\right]
	\\
	&\leq
	b' \norm{\n} a ^{\norm{\n}-1}
	\{
	c(\norm{\n}-1) \left[ (\norm{\n}-1)^2- a \norm{\n} ^2\right]
	\\
	&
	\quad\quad\quad\quad\quad\quad\quad\quad
	+a d \left[
	a (\norm{\n}+1)^2 + a^2 (\norm{\n}+2)^2 -2 \norm{\n}^2
	\right]
	\}
	\\
	& \leq b ,
	\end{aligned}
	\end{equation*}
	 where the last inequality holds since $a>1$, for some $b>0 $, and all $\n\in \mathbb{N}^M$.
The inequality in the last display, together with \eqref{UI}, implies
	\begin{equation}
	\label{UIfinal}
	\E{H(\N(t\wedge\tau_j))| \N(0)=\n}\leq 
	b' \norm{\n}^2 a ^{\norm{\n}} + b T.
	\end{equation}
Inequalities \eqref{fctH} and \eqref{UIfinal} finally ensure the method of duality can indeed be applied, which guarantees that the duality relationship between the generators \eqref{dualityrelationship}, proved in this section, implies the duality among the processes in the sense of \eqref{duality_expectations}. This completes the proof of Theorem \ref{thm:main}.
\\

\subsection{Proof of Corollary \ref{coroll}}
Assume $\theta_l=0$, $l=1,\dots,L$, and let $k$ be defined as in \eqref{k_nomut}. The rewriting of the diffusion and selection terms in \eqref{diffusion} and \eqref{selection} remains valid, even if \eqref{assumption0} is not satisfied in this case. 
Furthermore, it is possible to explicitly calculate,
for $l,r=1,\dots L, r\neq l, i,j=1,\dots,M_l, h=1,\dots,M_r$,
\begin{align*}
     \frac{k(\n-\e_i^{(l)})}{k(\n)}=\frac{n^{(l)}-1}{n_i^{(l)}-1}, \quad
    \frac{k(\n+\e_j^{(l)})}{k(\n)}  =\frac{n_j^{(l)}}{n^{(l)}}, \quad 
    \frac{k(\n+\e_j^{(l)}+\e_h^{(r)})}{k(\n)}  =\frac{n_j^{(l)}n_h^{(r)}}{n^{(l)}n^{(r)}} .
\end{align*}

Replacing these ratios in \eqref{diffusion} and \eqref{selection} yields an expression of the form \eqref{LF2} and provides the expression for the transition rates. As outlined in Section \ref{outline} an expression of the form \eqref{LF2} implies that the duality relationship between the generators of the diffusion and its dual process holds if 
\eqref{sumq} is satisfied. Since a stationary distribution that satisfies \eqref{assumption0} does not exist, the argument in Section \ref{outline} for proving \eqref{sumq} cannot be applied. However, direct calculation shows that
\begin{equation*}
\begin{aligned}
    &\sum_{l=1}^{L}\sum_{i=1}^{M_l} q(\n,\n-\e_i^{(l)})
    + \sum_{l=1}^{L}\sum_{j=1}^{M_l} q(\n,\n+\e_j^{(l)})
    + \sum_{l=1}^{L}\sum_{\substack{r=1\\r>l}}^{L}\sum_{j=1}^{M_l}\sum_{h=1}^{M_r}
    q(\n,\n+\e_j^{(l)}+\e_h^{(r)})
    \\
    &=
    \sum_{l=1}^{L}  
		\left(
		\frac{1}{2}n^{(l)}(n^{(l)}-1) 
		+
		\sum_{i=1}^{M_l}\sum_{\substack{k=1\\k\neq i}}^{M_l}
   	 	n_{i}^{(l)}h_{k}^{(l)}
    		+\sum_{\substack{r=1\\r\neq l}}^{L}\sum_{k=1}^{M_l}\sum_{m=1}^{M_r}	
    		n^{(l)}J_{km}^{(lr)}
		\right),
    \end{aligned}
    \end{equation*}
which implies \eqref{sumq}. 
Finally, the method of duality, using Corollary 4.4.3 in \cite{Ethier1986} as in the proof of Theorem \ref{thm:main}, ensures the duality relationship between the processes holds.

\section{Two loci, two alleles, with pairwise selection and parent independent mutation} 
\label{ex}
In this section a particular example will be considered, where there are two loci, $L=2$,  and two allele types at each locus, $M_1=M_2=2$. The pairwise interactions are represented by the matrix
\begin{align*}
J=
\left(\begin{array}{llll} 
0 & 0 &  J_1 & 0 \\
0 & 0 & 0 & J_2 \\
J_1 & 0 & 0 & 0 \\
0 & J_2 & 0 & 0 \end{array} \right), 
\end{align*}
and there is no single-locus selection,  $\mathbf{h}=0$. Furthermore, parent independent mutations are assumed.  

In this special case,  the function $k$,  in \eqref{eq:k}, and consequently the transition rates of the dual process can be computed rather efficiently.  The main difficulty in the computation is that the normalising constant of the stationary density \eqref{eq:stationary} is unknown. It 
may be noted that computing the normalising constant and the function $k$ are closely related problems. 
In fact, by defining
    \begin{equation*}
    I(a_1,a_2,b_1,b_2)=
    \int_0^1 \!\! \int_0^1 \!\!\!\!
    x^{a_1-1}(1-x)^{a_2-1} y^{b_1-1}(1-y)^{b_2-1}
    e^{ 2[J_1 xy  + J_2  (1-x) (1-y) ]}
    dx dy,
    \end{equation*}
the normalising constant can be written   as
    $$
    Z=I(2u_1^{(1)},2u_2^{(1)},2u_1^{(2)},2u_2^{(2)}), 
    $$
    and the function $k$ as
    $$
    k(\n)= \frac{1}{Z}
    I(n_1^{(1)}+2u_1^{(1)},n_2^{(1)}+2u_2^{(1)},n_1^{(2)}+2u_1^{(2)},n_2^{(2)}+2u_2^{(2)}).
    $$	
The integral $I$ cannot be computed analytically, but it is possible to find a series representation of it  in terms of Beta and Kummer functions, which can be  truncated to numerically evaluate the function $k$. The following formula is derived by a straightforward, albeit cumbersome, application of definitions and properties of Kummer functions
    \begin{align*}
    I(a_1,a_2 ,b_1,b_2) &=  e^{2J_2} B(a_1,a_2) 
    \sum_{n=0}^{\infty}  \frac{[a_1]_n}{[a_1+a_2]_n} \frac{(-2J_2)^n}{n!} \\
    & \quad \times \sum_{k=0}^n \binom{n}{k}\left(-\frac{J_1+J_2}{J_2}\right)^k
    B(k+b_1,b_2) {}_1F_1(k+b_1,k+b_1+b_2,-2J_2),
    \end{align*}
where $B$ is the Beta function, ${}_1F_1$ is the Kummer function and $[a]_n=a(a+1)\cdots (a+n-1)$, for $n>0$, $[a]_0=1$.

\begin{figure}[H]
	\centering
	\includegraphics[height=0.35\textwidth,width=\textwidth]{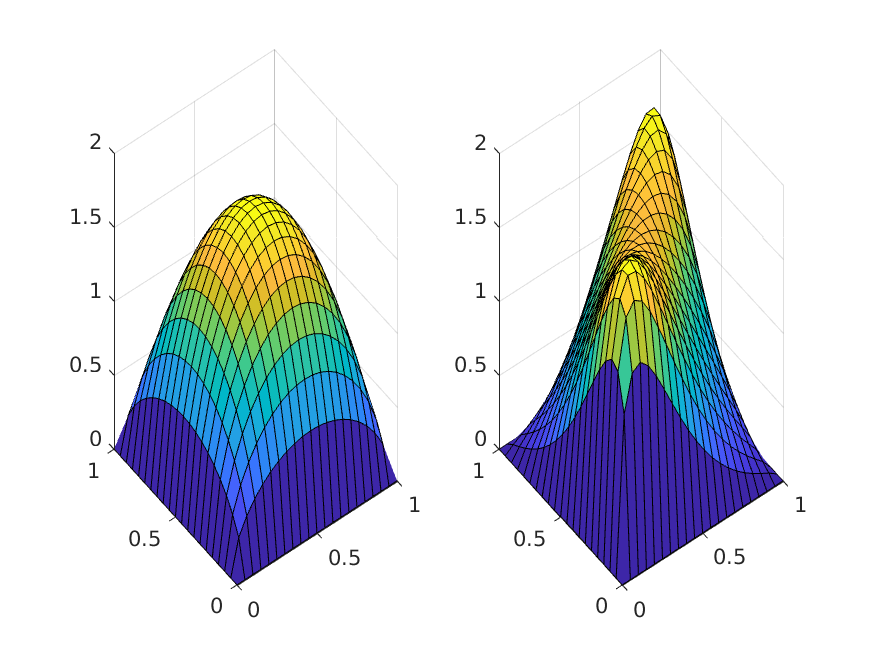}
	\caption{Stationary density of a  coupled Wright-Fisher diffusion for two loci, two alleles,  with no interaction (left)  
		and non-zero interaction (right). Mutation  parameters:
		$u_1^{(1)}=u_2^{(1)}=u_1^{(2)}=u_2^{(2)}=0.8$.
		Double-locus selection parameters:
		$J_1=J_2=0$ (left), $J_1=J_2=2$ (right).
	}
	\label{fig:statdens}
\end{figure}

As an illustration, the stationary density of independent Wright-Fisher diffusions, with $J_1 = J_2 = 0$,  is compared to the stationary density of the coupled Wright-Fisher diffusion, with $J_1 = J_2 = 2$,  in Figure \ref{fig:statdens}. Both distributions have mutation rates $u_1^{(1)}=u_2^{(1)}=u_1^{(2)}=u_2^{(2)}=0.8$. 
On the left hand side  the mutation strength keeps the mass of the stationary distribution in the centre of the unit square. In contrast, on the right hand side, while the mutation strength still tends to keep the mass in the centre, the selection strength moves the mass towards the points $(0,0)$ and $(1,1)$, which represent the most  viable couples of allele types, i.e.,  $1,1$ and $2,2$.

\section{Appendix}
In this section the monotone coupling mentioned in Section \ref{multi}  is presented.
Let $\{C_k\}_{k\in \mathbb{N}}\subset \mathbb{N}\setminus\{0\}$ be the jump chain of the process $C(t)=\norm{\N(t)}$, which represents the evolution of the number of genes in the genealogical (dual) process. 
Let $\{Y_k\}_{k\in \mathbb{N}} \subset \mathbb{N}$ be the Markov chain with the  following transition probabilities
	$$\mathbb{P}(Y_{k+1}=Y_k -1|Y_k)= p(Y_k), \quad \mathbb{P}(Y_{k+1}=Y_k +2|Y_k)= 1-p(Y_k), 
	$$
where $p(m)=\frac{c(m-1)}{c(m-1)+d}$, if $m>1$, $c,d$ are positive constants defined below in \eqref{ratesbounds}, $p(1)=p(2)$ and 
$p(0)=0.$
\\
As discussed in Section \ref{multi}, and shown in some particular cases in Sections \ref{ASG} and \ref{ex}, the rates of the dual process are often not explicit, furthermore, the process $C$ is not Markovian because its transition probabilities depend on $\N$. 
The aim of this appendix is to construct a coupling between  $\{C_k\}_{k\in \mathbb{N}}$ and $\{Y_k\}_{k\in \mathbb{N}}$ such that
	$$
	C_k \leq Y_k +1, \quad  \forall k\in \mathbb{N}.
	$$
This monotone coupling  provides upper bounds for expectations involving the number of genes in the genealogical process that is dual to the coupled Wright-Fisher diffusion. The Markov chain $Y$, with explicit and simple transition probabilities, is  easier to work with than  $C$. 

Let $\{\N_k\}_{k\in \mathbb{N}}$ be the jump chain of the process  $\N$, in which only coalescence and branching jumps are considered,
in between these jumps  the state of the process  $\N$ changes, because of an arbitrary number of mutations, from $\n$ to $\m$, with $\norm{\m}=\norm{\n}$,
say with probability $p(\m|\n)$. 
More precisely, given $C_k$ and $\N_k$,
$C_{k+1}$  is equal to
$C_k+j $ with probability $p_{j}(\N_k)$, $j=-1,1,2$, where

	\begin{equation}
	\label{transprob}
	p_{j}(\n)=\sum_{\m:\norm{\m}=\norm{\n}} p_{j}'(\m)p(\m|\n), 
	\quad \text{with}\quad
	p_{j}'(\m)=\frac{r_{j}(\m)}{r_{-1}(\m)+r_{1}(\m)+r_{2}(\m)} ,
	\end{equation}
and
	\begin{equation*}
	\begin{aligned}
	 r_{-1}(\n)&=\sum_{l=1}^L \sum_{\substack{i=1\\n_i^{(l)}>1}}^{M_l} q(\n,\n-\e_i^{(l)})	,\\
	 r_{1}(\n)&=\sum_{l=1}^L \sum_{\substack{j=1}}^{M_l} q(\n,\n+\e_j^{(l)}) ,
	\\
	 r_{2}(\n)&=\sum_{\substack{l,r=1\\l\neq r}}^L \sum_{\substack{j=1}}^{M_l} \sum_{\substack{j=1}}^{M_r}
	q(\n,\n+\e_j^{(l)}+\e_h^{(r)})	 .
	\end{aligned}
	\end{equation*}
The next step is to  bound the transition probabilities of  the chain $C$.
First note that definition \eqref{eq:k} yields
	$$
	k(\n - \e_i^{(l)}) \geq k(\n), 
	\quad 
	k(\n + \e_j^{(l)})\leq k(\n),	
	\quad
	k(\n + \e_j^{(l)}+\e_h^{(r)})\leq k(\n).$$
It is then straightforward to show that
	\begin{equation}
	\label{ratesbounds}
	r_{-1}(\n)
	\geq
	c \norm{\n}(\norm{\n}-1) ,
	\quad
	r_{1}(\n) +
	r_{2}(\n)
	\leq
	d \norm{\n},
	\end{equation}
where $c=\frac{1}{2\sqrt{M}}, d= 3\norm{J} +\norm{\mathbf{h}}$,
and thus that
	$$p_{-1}'(\n)\geq p(\norm{\n}), \quad p_1'(\n)+p_2'(\n) \leq 1- p(\norm{\n}).$$
This, together with \eqref{transprob}, implies
 	\begin{equation}
 	\label{probineq}
 	p_{-1}(\n)\geq p(\norm{\n}), \quad p_1(\n)+p_2(\n) \leq 1- p(\norm{\n}).
 	\end{equation}	
The inequalities above explain why  it is possible to construct a monotone coupling of $C$ and $Y$: 
$C$ has a higher probability of coalescence jump than $Y$ and a lower probability of branching jumps. Furthermore, inequalities \eqref{ratesbounds} are used in the proof of Theorem \ref{thm:main} to provide bounds that justify the use of the method of duality.

The coupling is  constructed step by step by coupling $C_k$ and $Y_k$ for each $k$ depending on which one  of the following cases occurs.

\textit{Case I}: $C_k=Y_k$.
In this case $C_{k+1}$ and $Y_{k+1}$ are coupled as follows, let $U$ be a standard uniform random variable, and set
	\begin{equation}
	\label{coupling}
	\begin{aligned}
	C_{k+1}=
	\begin{cases}
	C_k -1 \text{ if } U\in [0, p_{-1}(\N_k)], \\
	C_k+1 \text{ if } U\in (p_{-1}(\N_k),p_{-1}(\N_k)+p_{1}(\N_k)], \\
	C_k+2 \text{ if } U\in (p_{-1}(\N_k)+p_{1}(\N_k),1], \\
	\end{cases} 
	Y_{k+1}=
	\begin{cases}
	Y_k -1 \text{ if } U\in [0, p(Y_k)], \\
	Y_k+2 \text{ if } U\in (p(Y_k),1]. \\
	\end{cases}
	\end{aligned}
	\end{equation}
It is clear that, in this construction, the marginal distributions are preserved and additionally  $C_{k+1}\leq Y_{k+1}$ because of \eqref{probineq} and $\norm{\N_k}=C_k=Y_k.$

\textit{Case II}:
$C_k<Y_k$. As long as  this case occurs,  let $C$ and $Y$ evolve independently.

\textit{Case III}:
$C_k>Y_k$. Assume that $k$ is first time this case occurs again after case II (it cannot occur after case I) and note that $C_k=Y_k+1$ must hold. In this case $C_{k+1} $ and $Y_{k+1}$ are coupled  as in \eqref{coupling}, the difference being that here $\norm{\N_k}=C_k=Y_k +1$. Since  $p$ is an increasing function,  $C_{k+1}\leq Y_{k+1}+1$ holds. This means that after  one step in case III, either case I occurs, or  $C_{k+1}= Y_{k+1}+1$ and the latter coupling  can be applied again. 

Note that coupling $C_k$ and $Y_k$   impose implicitly a coupling on $\N_k$ and $Y_k$.
Applying the appropriate coupling at each step provides a coupling between the chain $C$ and the Markov chain $Y$ such that $C_{k}\leq Y_{k}+1$, $\forall k \in \mathbb{N}$, assuming that $C_0=Y_0=\norm{\N(0)}$. 
Furthermore, it is interesting to note that the first time $C$ hits $1$, which is the time the genealogical process reaches the most recent common ancestor,  is smaller or equal to the first time $Y$ hits $0$. In fact,
 either $C$ hits $1$ before $Y$ hits $0$, or when $Y$ hits $0$ it must be that $C $ hits $1$.

\begin{acknowledgements}
We are very grateful to the associate editor and two anonymous referees whose detailed comments and suggestions lead to significant improvements of the manuscript. The research is supported by the Swedish Research Council through grants 621-2013-4628 (H. Hult) and 40-2012-5952 (T. Koski).  
\end{acknowledgements}

%
%

\bibliographystyle{spbasic}      
\bibliography{mybib}   

\end{document}